\documentclass{amsart}
\usepackage{graphicx}
\usepackage{float}
\usepackage[arc,all]{xy}
\usepackage{longtable}
\newtheorem{thm}{Theorem}[section]

\newtheorem{lem}[thm]{Lemma}
\newtheorem{prop}[thm]{Proposition}
\theoremstyle{definition}

\newtheorem*{proof M}{\textbf{Proof of  Main Theorem}}
\theoremstyle{remark}

\numberwithin{equation}{section}

\begin{document}
\title[the tensor product of nilpotent Lie algebras ]{On the triple tensor product of nilpotent Lie algebras }%
\author{afsaneh shamsaki}%
\address{School of Mathematics and Computer Science\\
Damghan University, Damghan, Iran}
\email{Shamsaki.afsaneh@yahoo.com}%
\author[P. Niroomand]{Peyman Niroomand}\thanks{Corresponding Author: Peyman Niroomand}
\email{niroomand@du.ac.ir, p$\_$niroomand@yahoo.com}
\address{School of Mathematics and Computer Science\\
Damghan University, Damghan, Iran}

\begin{abstract}
In this paper,  we give the explicit structure of $ \otimes^{3} H $ and  $ \wedge^{3} H $ where $ H $ is a  generalized Heisenberg Lie algebra of rank at most  $ 2. $ Moreover,  for a non-abelian nilpotent Lie algebra $ L, $ we obtain an upper bound for the dimension of $  \otimes^{3} L. $ 
\end{abstract}
\maketitle
\section{ Introduction}
It is known from \cite[Theorem 3.1]{N} that the dimension of tensor product of a Lie algebra $ L $ of dimension $ n $ when $ \dim L^{2}=m $ is bounded by $ (n-m)(n-1)+2 $ and the structure of the tensor product of a Heisenberg Lie algebra is given. Let $ H $ be a generalized Heisenberg Lie algebra of rank at most $ 2.  $
In this paper, we complete  the explicit structure of the tensor square and the exterior square of $ H. $
Moreover, we give the triple tensor product and the triple exterior product of $H. $ Finally, we show that dimension of the triple tensor of a Lie algebra $ L $ with derived subalgebra of dimension $ m $ is bounded by $ n(n-m)^{2}. $\\
\section{ Preliminary and known results}
Let $ \mathbb{F} $ be a fixed field and $ [, ] $ is used to denote the Lie bracket. For any two Lie algebras $ L $ and $ K, $ we say there exists an action $ L $ on $ K $ if an $ \mathbb{F} $-bilinear map  $ L\times K\rightarrow K, $ $ (l, k)\mapsto ~^{l}k $ satisfying 
\begin{align*}
^{[l, l']}k= ^{l}(^{l'}k)-^{l'}(^{l}k)\qquad \text{and}\qquad ^{l}[k, k']=[^{l}k, k']+[k, ^{l}k']
\end{align*}
for all $ l, l'\in L $ and $ k, k' \in K. $ \\
The actions are compatible if
\begin{equation}\label{e}
^{^{l}k}l'=[l', ^{k}l] \qquad \text{and} \qquad ^{^{k}l}k'=[k', ^{l}k]
\end{equation}
for all $ k, k' \in K, $ $ l,l'\in L. $ \\
Clearly, a Lie algebra $ L $ acts on itself by   $ ^{l'}l=[l',l] $ for all $ l, l' \in L. $ Morover,  for all Lie algebras $ Q, $
a bilinear function $ \varphi : L\times K\rightarrow Q $ is called a Lie pairing if for all $ l, l'\in L $ and $ k, k'\in K, $ we have
\begin{align*}
& \varphi ([l, l'], k)=\varphi (l, ^{l'}k)-\varphi(l', ^{l}k),\cr
& \varphi(l, [k, k'])=\varphi(^{k'}l,k)-\varphi (^{k}l, k'),\cr
& \varphi (^{k}l, ^{l'}k')=-[\varphi(l,k), \varphi(l',k')].
\end{align*} 
The non-abelian tensor product $ L\otimes K $ is the Lie algebra generated by symbols $ l\otimes k $ for all $ l\in L $ and $ k\in K $ with the following defining relations
\begin{align*}
&c(l\otimes k)=cl\otimes k=l\otimes ck, \cr
& (l+l')\otimes k=l\otimes k+l'\otimes k,\cr
& l\otimes (k+k')=l\otimes k+l\otimes k',\cr
& [l, l']\otimes k=l\otimes ^{l'}k-l'\otimes~ ^{l}k,\cr
& l\otimes [k, k']= ^{k'}l\otimes k-^{k}l\otimes k',\cr
& [(l\otimes k), (l\otimes k')]=-^{k}l\otimes ~^{l'}k'
\end{align*} 
for all $ c\in \mathbb{F}, $ $ l, l'\in L $ and $ k, k' \in K. $ If $ L=K $ and all actions are Lie multiplication, then $ L\otimes L $ is called the non-abelian tensor square of $ L. $\\
The next two propositions are used in the rest of paper.
\begin{prop}\cite[Proposition 3]{E}\label{p1.1}
There are actions of both $ L $ and $ K $ on $L\otimes K $ given by
\begin{align*}
&^{l'}(l\otimes k)=[l, l']\otimes k+l\otimes (^{l'}k),\cr
& ^{k'}(l\otimes k)=(^{k'}l)\otimes k+l\otimes [k', k]
\end{align*}  
for all $ l, l' \in L $ and $ k, k' \in K. $
\end{prop}
Recall that a Lie algebra $ H $   is called a generalized Heisenberg of rank $ n $ if $ H^{2}=Z(H) $ and $ \dim H^{2}=n. $ If $ n=1, $ then $ H $ is a Heisenberg Lie algebra that is more well-known. Such Lie algebras are odd dimension and $ H\cong H(m)= \langle x_{i}, y_{i}, z\mid [x_{i}, y_{i}]=z, 1\leq i \leq m \rangle.$
\begin{prop}\cite[Proposition 3.3]{N}\label{p1.2}
Let $ H(m) $ be a Heisenberg Lie algebra. Then 
\begin{equation*}
H(m)\otimes H(m)\cong H(m)/H(m)^{2}\otimes H(m)/H(m)^{2}
\end{equation*} 
for all $ m $ such that $ m\geq 2. $ In the case $ m=1, $ $ H(1)\otimes H(1) $ is an abelian Lie algebra of dimension $ 6. $ 
\end{prop}
Let $ L\square L $ be the submodule of $ L\otimes L $ generated by the elements $ l\otimes l. $ Then the exterior square $ L\wedge L $ of $ L $ is the quotient $ L\wedge L\cong L\otimes L/ L\square L. $ For all $ l\otimes l'  \in L\otimes L,$ we denote the coset $ l\otimes l'+L\square L $ by $ l\wedge l'. $ \\
The Schur multiplier of a Lie algebra $ L $ is defined as $ \mathcal{M}(L)\cong R\cap F^{2}/[R, F] $ where $ L\cong F/R $ and $ F $ is a free Lie algebra. It is known  that the Lie algebra $ \mathcal{M}(L) $ is abelian and independent of the choice of the free Lie algebra $ F$ (see \cite{B1, B2} for more information).
\begin{lem}\cite[Lemma 23]{M}\label{M}
Let $ L $ be an $ n $-dimensional abelian Lie algebra. Then $ \dim \mathcal{M}(L)=\frac{1}{2}n(n-1).  $
\end{lem}
\begin{lem}\cite[Lemma 2.8 ]{J}\label{J}
Let $ L $ be a nilpotent Lie algebra of class two. Then $ L\wedge L\cong \mathcal{M}(L)\oplus L^{2}. $
\end{lem}
Let $ A(n) $ be used to denote an abelian Lie algebra of dimension $ n. $ Then we have
\begin{lem}\cite[Corollary 2.5 and Lemma 3.2]{N3}\label{H}
\begin{itemize}
\item[(i).] $ A(n)\wedge A(n)\cong \mathcal{M}(A(n)),  $
\item[(ii).] $ H(1)\wedge H(1)\cong A(3), $
\item[(iii).] $ H(m)\wedge H(m)\cong A(2m^{2}-m) $ for all $ m\geq 2. $
\end{itemize}
\end{lem}

The next theorem shows that the exterior product of a Lie algebra $ L $ is a direct summand of $ L\otimes L $ when $ L/L^{2} $ is of finite dimension.
\begin{thm}\cite[Theorem 2.5]{N1}\label{1.5}
Let $ L/L^{2} $ be a finite dimensional Lie algebra. Then
\begin{equation*}
L\otimes L\cong L\wedge L\oplus L\square L.
\end{equation*}
\end{thm}
In the following, the notations $ L_{6,22}(\varepsilon), $      $L_{5,8}, $  $L_{6, 7}^{(2)}(\eta)$ and $ L_{1} $  are  taken from \cite{G} and \cite{N2}.
\begin{prop}\cite[Theorem 2.15]{J}
The Schur multiplier of Lie algebras $ L_{6,22}(\varepsilon), $  $ L_{5,8}, $ $L_{6, 7}^{(2)}(\eta),$   and $ L_{1} $ are abelian Lie algebras of dimension $ 8, 6, 8 $ and $ 9, $ respectively.  
\end{prop}
A Lie algebra $ L $ is called capable provided that $ L\cong H/Z(H) $ for a Lie algebra $ H. $ 
\begin{prop}\cite[Proposition 3.1]{N2}\label{Shur2}
Let $ H $ be a non-capable generalized Heisenberg Lie algebra of rank $ 2 $ such that $ \dim H=n. $ Then
\begin{align*}
\dim \mathcal{M}(H)=\frac{1}{2}(n-3)(n-2)-2
\end{align*}
or
\begin{align*}
\dim \mathcal{M}(H)=\frac{1}{2}(n-1)(n-4)+1.
\end{align*}
\end{prop}
\section{main results}
We know from Proposition   \ref{p1.1}, $ L $ acts on $ L\otimes L. $ On the other hand,  the tensor product $ L\otimes L $ acts on $ L $ by $ ^{t}l=~~ ^{\lambda (t)}l$ for all $ t\in L\otimes L $ and $ l\in L $ such that $ \lambda : L\otimes L\rightarrow L $ is a homomorphism  given by   $ a\otimes b\mapsto [a, b].$ These actions are compatible and we can construct the triple tensor product $ \otimes^{3}L=(L\otimes L)\otimes L. $ \\
In this section, we give the explicit structure of $ L\otimes L=\otimes^{2} L $ and  $ L\wedge L=\wedge^{2} L$ when $ L $ is a generalized Heisenberg of rank $ 2. $ Moreover, we obtain $  \otimes^{3}L $ and $ \wedge^{3}L $ when $ L $ is a generalized Heisenberg of rank at most $ 2 $. Also, for a non-abelian nilpotent Lie algebra $ L, $ we  give an upper bound for the triple tensor product of $ L. $ 
\\The following lemmas are  useful instruments in the next.
\begin{lem}\label{l1}
Let $ L $ be a Lie algebra of nilpotency class two. Then
\begin{itemize}
\item[(i).] $ L\otimes L $ acts trivially on $ L. $
\item[(ii).] $ (L\otimes L)\otimes L $ is an abelian Lie algebra.
\end{itemize}
\begin{proof}
(i). We know that $ L $ is a Lie algebra of nilpotency class two, so $ \gamma_{3}(L)=0. $ By considering action of   $ L\otimes L $  on $ L, $ we have $ ^{\lambda(a\otimes b)}c=[[a, b], c]=0 $ for all $ a\otimes b\in L\otimes L $ and $ c\in L. $ Therefore $ L\otimes L $ acts trivially on $ L. $\\
(ii). By using the relations of the non-abelian tensor on $ (L\otimes L)\otimes L $ and the part (i), we have
\begin{align*}
[(a\otimes b)\otimes c, (a'\otimes b')\otimes c']&= -(^{c}(a\otimes b)) \otimes ^{a'\otimes b'}c\cr
&=-^{c}(a\otimes b)\otimes 0=0
\end{align*} 
Therefore $ [(a\otimes b)\otimes c, (a'\otimes b')\otimes c']=0 $ for all $ (a\otimes b)\otimes c, (a'\otimes b')\otimes c' \in (L\otimes L)\otimes L, $ and so $ (L\otimes L)\otimes L $ is an abelian Lie algebra.
\end{proof}
\end{lem}
Let $ L $ be a nilpotent Lie algebra of class $ k $, $ i_{L}: L\rightarrow L $ be the identity homomorphism and  $ \gamma_{k}(L) $ be the $ k $-th term of the lower central series of $ L $ and $ \varphi :  \gamma_{k}(L) \rightarrow L$ be a natural homomorphism.  Define homomorphisms  
 $ \overline{\varphi}=(\varphi \otimes i_{L})\otimes i_{L}  :  (\gamma_{k}(L)\otimes L)\otimes L \rightarrow \otimes^{3}L $
 and   $ \gamma : (L\otimes L) \otimes \gamma_{k}(L) \rightarrow  \otimes^{3} L $ given by  $ (a\otimes b)\otimes c\mapsto (a\otimes b) \otimes c $. Then
 
\begin{lem}\label{l2}
Let $ L $ be a non-abelian nilpotent Lie algebra of class $ k, $ then  $ \mathrm{Im} \gamma \subseteq \mathrm{Im}  \overline{\varphi}. $
\begin{proof}
  By using the relation $  k\otimes [l, l']=[l', k] \otimes l- [l, k] \otimes l'$
   on $ L\otimes L, $ we can see that  
$ (a\otimes b) \otimes [x_{1}, \dots, x_{k-1}, x_{k}]\in  (\gamma_{k}(L)\otimes L)\otimes L$
 for all $ a, b, x_{1}, \dots, x_{k-1}, x_{k} \in L.$ Therefore the result follows.
\end{proof}
\end{lem}
The next result plays a key role in proving the next theorem.
 \begin{prop}\label{p2}
 If $ L $ is a nilpotent Lie algebra of class $ k, $ then 
 \begin{equation*}\label{ee}
 (\gamma_{k}(L)\otimes L)\otimes L  \overset{(\varphi \otimes i_{L})\otimes i_{L}}\longrightarrow \otimes^{3}L \longrightarrow \otimes^{3} L/\gamma_{k}(L)\rightarrow 0,
 \end{equation*}
 is exact.
 \begin{proof}
By using  Lemma \ref{l2} and  \cite[Proposition 1.3]{sal}, the result follows. 
 \end{proof}
 \end{prop}

 In the following proposition, we are going to determine the structure of $ \otimes^{3} H(m) $ and $ \wedge^{3} H(m) $  for all $ m\geq 1. $ 
\begin{prop}\label{p2.3}
Let $ H(m) $ be a Heisenberg Lie algebra. Then 
$$
\otimes^{3}H(m)\cong 
\begin{cases}
  A(12) &\text{if }m=1,\\
 A(2^{3}m^{3}) &\text{if } m\geq2\\
\end{cases}
$$
and
$$
\wedge^{3}H(m)\cong 
\begin{cases}
  A(2) &\text{if }m=1,\\
 A(4m^{3}-4m^{2}-m) &\text{if } m\geq2.\\
\end{cases}
$$
 \begin{proof}
We claim that $( \varphi(L^{2})\otimes L)\otimes L \cong  (\varphi(L^{2})\otimes L^{ab})\otimes L^{ab}.$ It is clear that $ \varphi(L^{2})$ and $ L $ act trivially on each other. Hence  $( \varphi(L^{2})\otimes L)\otimes L \cong  (\varphi(L^{2})\otimes L^{ab})\otimes L.$ Also, $ \varphi(L^{2})\otimes L^{ab} $ and $ L $ act trivially on each other. Thus $( \varphi(L^{2})\otimes L)\otimes L \cong  (\varphi(L^{2})\otimes L^{ab})\otimes L^{ab}.$ By using the following exact sequence
\begin{equation*}
L^{2}\otimes L\overset{\varphi \otimes i_{L}}\longrightarrow L\otimes L\longrightarrow L/L^{2}\otimes L/L^{2}\longrightarrow 0,
\end{equation*}
 we have  
\begin{equation}\label{e3}
\dim L\otimes L=\dim L^{ab}\otimes L^{ab}+\dim \mathrm{Im} \varphi \otimes i_{L} =\dim L^{ab}\otimes L^{ab}+\dim \varphi(L^{2})\otimes L.
\end{equation}
Let $ L\cong H(1). $ Since $ H(1)\otimes H(1)\cong A(6) $ and $ H(1)^{ab}\otimes H(1)^{ab}\cong A(4), $ we have $ \dim \varphi(L^{2})\otimes L=2 $ by using (\ref{e3}). Since  $( \varphi(L^{2})\otimes L)\otimes L \cong  (\varphi(L^{2})\otimes L^{ab})\otimes L^{ab}$,
\begin{equation}\label{e4}
\dim  (\varphi(L^{2})\otimes L)\otimes L =  (\dim \varphi(L^{2})\otimes L^{ab})\dim L^{ab}.
\end{equation}
Now, by using (\ref{e4}), Proposition \ref{p2} and  Lemma \ref{l1} (ii), we have
\begin{equation*}
\dim \otimes^{3} L=\dim \otimes^{3} L^{ab}+ \dim (\varphi(L^{2})\otimes L^{ab})\otimes L^{ab}=12
\end{equation*} 
and so $ \otimes^{3} L\cong A(12). $ Similarly, if $ L\cong H(m) $ for all $ m\geq 2, $ then $ \otimes^{3}L\cong A(2^{3}m^{3}). $\\
Let $ L\cong H(1), $   then $ \dim \varphi(L^{2})\wedge L= \dim\wedge^{2} L-\dim \wedge^{2} L^{ab}=  2$ by using the following exact sequence 
\begin{equation*}
L^{2}\wedge L\overset{\varphi \wedge i_{L}}\longrightarrow L\wedge L\longrightarrow L/L^{2}\wedge L/L^{2}\longrightarrow 0
\end{equation*}
 and Lemma \ref{H}. We can easily see $ (\varphi(L^{2})\wedge L) \wedge L\cong (\varphi(L^{2})\wedge L^{ab}) \wedge L^{ab} $ and since $ \varphi(L^{2})\wedge L\cong A(2),  $ we have $ (\varphi(L^{2})\wedge L) \wedge L\cong A(2)\wedge A(2). $ Thus $ (\varphi(L^{2})\wedge L) \wedge L\cong A(1) $ by using  Lemma \ref{H}. On the other hand, $ \wedge^{3} L/L^{2}\cong (A(2)\wedge A(2))\wedge A(2), $ hence $ \wedge^{3} L/L^{2}\cong A(1) $ by using Lemma \ref{H} and \cite[Lemma 2.1.7]{J1}.
Consider the following exact sequence
\begin{equation*}
 (L^{2} \wedge L)\wedge L  \overset{(\varphi\wedge i_{L})\wedge i_{L}}\longrightarrow \wedge^{3}L \longrightarrow \wedge^{3} L/L^{2}\rightarrow 0.
\end{equation*}
Therefore
\begin{align*}
 \dim \wedge^{3}L &=\dim\wedge^{3} L/L^{2}+\dim \mathrm{Im} (\varphi \wedge i_{L})\cr
&=\dim \dim\wedge^{3} L/L^{2}+\dim (\varphi(L^{2})\wedge L) \wedge L\cr
&=2.
\end{align*}
Since $ \otimes^{3} L $ is an abelian Lie algebra by using Lemma \ref{l1} (ii), we have $ \wedge^{3} L\cong A(2). $ Similarly, if $ L\cong H(m) $ for all $ m\geq 2, $ then $ \wedge^{3}L\cong A(4m^{3}-4m^{2}-m). $\\
 \end{proof}
\end{prop}
The following theorem determines  all $ n $-dimensional capable generalized Heisenberg Lie algebras of rank $ 2. $ 
\begin{thm}\cite[Theorem 3.6]{N2}
Let $ H $ be a generalized Heisenberg Lie algebra such that  $ \dim H=n $ and $ \dim H^{2}=2. $ Then $ H $ is capable if and only if $ n=5, 6 $ or $ 7 $ and $ H\cong L_{5,8}=\langle  x_1,\ldots, x_5\big{|}[x_1, x_2] = x_4, [x_1, x_3] = x_5\rangle, $ $ H\cong L_{6,22}(\varepsilon)=\langle   x_1,\ldots, x_6\big{|}[x_1, x_2] = x_5= [x_3, x_4], [x_1, x_3] = x_6,[x_2, x_4] = \epsilon x_6 \rangle $ where $ \varepsilon \in \mathbb{F}/(\stackrel{*}{\sim}), $  $ L_{6, 7}^{(2)}(\eta)=\langle x_1,\ldots, x_6 \mid [x_{1}, x_{2}]=x_{5},  [x_{1}, x_{3}]=x_{6}, [x_{2}, x_{4}]=\eta x_{6}, [x_{3}, x_{4}]=x_{5}+x_{6}$ where $ \eta \in \lbrace 0, \omega \rbrace  $
or $ H\cong L_{1}=\langle x_1,\ldots, x_7\big{|}[x_1, x_2] = x_6=[x_3, x_4],[x_1, x_5] = x_7= [x_2, x_3]\rangle. $  
\end{thm}
 The next  propositions  determine the structure of $ \otimes^{2} H, $  $ \wedge^{2}H, $ $ \otimes^{3} H $ and $ \wedge^{3}H $ when  $ H $ is a generalized Heisenberg Lie algebra of rank $ 2. $

 \begin{prop}\label{p'1}
Let $  H$ be an $ n $-dimensional generalized Heisenberg Lie algebra of rank $ 2. $
\begin{itemize}
\item[(i).] If $ H$ is non-capable, then 
\end{itemize}
$$
\otimes^{2}H \cong 
\begin{cases}
  A((n-2)^{2}) &\text{if }Z^{\wedge}(H)=H^{2},\\
 A((n-2)^{2}+2) &\text{otherwise} \\
\end{cases}
$$
and 
$$
\wedge^{2}H \cong 
\begin{cases}
 A(\frac{1}{2}(n-3)(n-2)) &\text{if }Z^{\wedge}(H)=H^{2},\\
A(\frac{1}{2}(n-1)(n-4)+3) &\text{otherwise.} \\
\end{cases}
$$
\item[(ii).] If $ H $ is capable, then
$$
\otimes^{2}H\cong 
\begin{cases}
  A(14) &\text{if }H\cong L_{5,8},\\
 A(20) &\text{if } H\cong L_{6,22}(\varepsilon)~~ \text{or}~~L_{6, 7}^{(2)}(\eta),\\
 A(26) &\text{if } H\cong L_{1}\\
\end{cases}
$$
and
$$
\wedge^{2}H\cong 
\begin{cases}
  A(8) &\text{if }H\cong L_{5,8},\\
 A(10) &\text{if } H\cong L_{6,22}(\varepsilon)~~ \text{or}~~L_{6, 7}^{(2)}(\eta),\\
 A(11) &\text{if } H\cong L_{1}.\\
\end{cases}
$$
\begin{proof}
(i). First, we obtain the exterior product of $ H. $   By invoking Lemma \ref{J} and Proposition \ref{Shur2}, we have $ \wedge^{2}H=\mathcal{M}(H)\oplus H^{2}\cong   A(\frac{1}{2}(n-3)(n-2))$ when $ Z^{\wedge}(H)=H^{2}. $  By a similar way, we can see that  $\wedge^{2}H= A(\frac{1}{2}(n-1)(n-4)+3) $ when $ Z^{\wedge}(H)\neq H^{2}. $ \\
Since $ H\square H\cong  H^{ab}\square H^{ab}\cong A(\frac{1}{2}(n-m)(n-m+1))$ by \cite[Lemma 2.3]{N1} and  using  Theorem \ref{1.5}, we have  $ \otimes^{2}H=H\square H\oplus H\wedge H \cong   A(n^{2}-2n+3) $ when  $ Z^{\wedge}(H)=H^{2}. $   By a similar way, we can see that $  \otimes^{2}H \cong A(n^{2}-2n+5) $ when $ Z^{\wedge}(H)\neq H^{2}. $ \\
(ii). The proof is  similar to  the part (i).
\end{proof}
\end{prop}
\begin{prop}\label{p1}
Let $  H$ be $ n $-dimensional generalized Heisenberg Lie algebra of rank $ 2. $
\begin{itemize}
\item[(i).] Let $H $ is non-capable. Then 
\end{itemize}
$$
\otimes^{3}H \cong 
\begin{cases}
  A((n-2)^{3}) &\text{if }Z^{\wedge}(H)=H^{2},\\
 A((n^{2}-4n+6)(n-2)) &\text{otherwise.} \\
\end{cases}
$$
and
$$
\wedge^{3}H \cong 
\begin{cases}
 A(\frac{1}{2}(n-2)(n^{2}-6n+7)) &\text{if }Z^{\wedge}(H)=H^{2},\\
A(\frac{1}{2}n(n^{2}-8n+23)-14)&\text{otherwise.} \\
\end{cases}
$$
\item[(ii).] If $H $ is capable, then
$$
\otimes^{3}H\cong 
\begin{cases}
  A(42) &\text{if }H\cong L_{5,8},\\
 A(80) &\text{if } H\cong L_{6,22}(\varepsilon)~~ \text{or}~~L_{6, 7}^{(2)}(\eta),\\
 A(130) &\text{if } H\cong L_{1}\\
\end{cases}
$$
and
$$
\wedge^{3}H\cong 
\begin{cases}
  A(12) &\text{if }H\cong L_{5,8},\\
 A(20) &\text{if } H\cong L_{6,22}(\varepsilon)~~ \text{or}~~L_{6, 7}^{(2)}(\eta),\\
 A(39) &\text{if } H\cong L_{1}.\\
\end{cases}
$$
\begin{proof}
The proof is obtained by a similar   way used  in the   proof of Proposition \ref{p2.3}.
\end{proof}
\end{prop}

 \begin{thm}
Let $ L $ be an $ n $-dimensional non-abelian nilpotent Lie algebra with derived subalgebra of dimension 
$ m. $ Then  
\begin{equation*}
\dim \otimes^{3} L  \leq n(n-m)^{2}.
\end{equation*} 
In particular, for $ m=1 $ the equality holds if and only if $ L\cong H(1). $
\begin{proof}
Let $ L $ be an $ n $-dimesional  nilpotent Lie algebra of nilpotency class $ k. $
We proceed by induction on $ n. $ Since $ L $ is non-abelian, $ n\geq 3. $  If $ n=3, $ then 
$ L\cong H(1) $ by the classification of all nilpotent Lie algebras is given  in \cite{G}. By using Proposition \ref{p2.3}, we have $ \otimes^{3} H(1)\cong A(12). $ Therefore the result holds. Since $ \gamma_{k}(L) $ is central,   $ \gamma_{k}(L) $ and $ L $ act trivially on each other and so we have $ \gamma_{k}(L)\otimes L\cong \gamma_{k}(L)\otimes L^{ab}$ by using \cite[Proposition 5]{E}. By the same reason and using \cite[Proposition 5]{E}, we have $ (\gamma_{k}(L)\otimes L)\otimes L\cong (\gamma_{k}(L)\otimes L^{ab})\otimes L^{ab}. $ 
   Consider two cases.\\
Let $ L/\gamma_{k}(L) $ is abelian. Then $ L $ is a  nilpotant Lie algebra of class two and  $ \gamma_{k}(L)=L^{2}. $ Thus $ \dim \otimes^{3} L/\gamma_{k}(L)=\dim \otimes^{3} L/L^{2}=(n-m)^{3}. $
By using  Proposition \ref{p2}, we have
\begin{align*}
\dim  \otimes ^{3}L & \leq \dim \otimes ^{3}(L/\gamma_{k}(L))+\dim ((\gamma_{k}(L)\otimes L)\otimes L)\cr
& =(n-m)^{3}+m(n-m)^{2}=n(n-m)^{2}.
\end{align*}
	If $ L/\gamma_{k}(L) $ is non-abelian, then $ \dim L/\gamma_{k}(L)=n-\dim \gamma_{k}(L) $ and $ \dim (L/\gamma_{k}(L))^{2}=m-\dim \gamma_{k}(L) , $ hence 
\begin{equation}\label{e2}
\dim  \otimes ^{3}(L/\gamma_{k}(L)) \leq (n-\dim \gamma_{k}(L))(n-m)^{2}
\end{equation}
 by using the induction hypothesis. 
Now, Proposition \ref{p2} and (\ref{e2}) imply that 
\begin{align*}
\dim \otimes ^{3}L&\leq \dim \otimes ^{3}(L/\gamma_{k}(L))+\dim ((\gamma_{k}(L)\otimes L)\otimes L)\cr
&\leq (n-\dim \gamma_{k}(L))(n-m)^{2}+\dim \gamma_{k}(L)(n-m)^{2} \cr
& = n(n-m)^{2}.
\end{align*}
Let $ m=1. $ Then $ L\cong H(k) \oplus A(n-2k-1)$ for all $ k\geq 1 $ by using \cite[Lemma 3.3]{N4}.
By using Proposition \ref{p2.3} the equality holds if and only if $ L\cong H(1). $ 
\end{proof}
 \end{thm}


\begin{thebibliography}{9}
\bibitem{B1} 
Batten, P., Moneyhun, K. and Stitzinger, E.;  On characterizing nilpotent Lie algebras by their multipliers. Comm. Algebra 24 (1996), no. 14, 4319--4330.
\bibitem{B2}
Batten, P. and  Stitzinger, E.; On covers of Lie algebras. Comm. Algebra 24 (1996), no. 14, 4301--4317.
\bibitem{E}
Ellis, G.;  A non-abelian tensor product of Lie algebras. Glasgow Mathematical Journal 33, no. 1 (1991), 101--120.
\bibitem{G}
De Graaf, W. A.;  Classification of 6-dimensional nilpotent Lie algebras over fields of charac-teristic not 2. J. Algebra 309, no. 2, (2007), 640–653.
\bibitem{M}
Moneyhun, K.; Isoclinisms in Lie algebras. Algebras Groups Geom. 11 , no. 1, (1994), 9--22.
\bibitem{J}
Johari, F., and  Niroomand, P.; Certain functors of nilpotent Lie algebras with the derived subalgebra of dimension two. J. Algebra Appl. 19, no. 01, (2020), 2050012.
\bibitem{J1}
Johari, F., Parvizi, M., and  Niroomand, P. Capability and Schur multiplier of a pair of Lie algebras. J. Geom. Phys. 114, (2017), 184--196.
\bibitem{N}
Niroomand, P.; On the tensor square of non-abelian nilpotent finite-dimensional Lie algebras. Linear and Multilinear Algebra 59, no. 8 (2011), 831--836.
\bibitem{N1}
Niroomand, P.; Some properties on the tensor square of Lie algebras. J. Algebra Appl. 11 (2012), no. 5, 1250085, 6 pp.
\bibitem{N2}
P. Niroomand, F. Johari and M. Parvizi, Capable Lie algebras with the derived subalgebra of dimension 2 over an arbitrary field. Linear and Multilinear Algebra 67 , no. 3, (2019), 542-554. 
 \bibitem{N3}
 Niroomand, P., Parvizi, M. and Russo, F. G.; Some criteria for detecting capable Lie algebras. J. Algebra 384 (2013), 36–44.
 \bibitem{N4}
 Niroomand, P.; On dimension of the Schur multiplier of nilpotent Lie algebras. Cent. Eur. J. Math. 9 (2011), no. 1, 57–64.
\bibitem{sal}
Salemkar, A.R., Tavallaee, H., Mohammadzadeh, H. and Edalatzadeh, B.; On the non-abelian tensor product of Lie algebras. Linear and Multilinear Algebra 58, no. 3, (2010), 333--341.

\end{thebibliography}
\end{document}